\documentclass[11pt]{article}
\usepackage{amsfonts,amssymb,amsmath}

\usepackage[english]{babel}
\usepackage{color}

\usepackage[pagebackref=true,colorlinks,allcolors=blue,breaklinks]{hyperref}

\newtheorem{theorem}{Theorem}
\newtheorem{definition}{Definition}
\newtheorem{remark}{Remark}
\newtheorem{proof}{Proof}

\newtheorem{proposition}{Proposition}
\newtheorem{lemma}{Lemma}

\newcommand{\beq}{\begin{eqnarray}}
\newcommand{\eeq}{\end{eqnarray}}

\newcommand{\beqt}{\begin{eqnarray*}}
\newcommand{\eeqt}{\end{eqnarray*}}

\newcommand{\be}{\begin{equation}}
\newcommand{\ee}{\end{equation}}

\newcommand{\bl}{\begin{lemma}}
\newcommand{\el}{\end{lemma}}

\newcommand{\bcon}{\begin{conjecture}}
\newcommand{\econ}{\end{conjecture}}

\newcommand{\br}{\begin{remark}}
\newcommand{\er}{\end{remark}}

\newcommand{\bt}{\begin{theorem}}
\newcommand{\et}{\end{theorem}}

\newcommand{\bd}{\begin{definition}}
\newcommand{\ed}{\end{definition}}

\newcommand{\bp}{\begin{proposition}}
\newcommand{\ep}{\end{proposition}}

\newcommand{\bc}{\begin{corollary}}
\newcommand{\ec}{\end{corollary}}

\newcommand{\bpr}{\begin{proof}}
\newcommand{\epr}{\end{proof}}

\newcommand{\bi}{\begin{itemize}}
\newcommand{\ei}{\end{itemize}}

\newcommand{\ben}{\begin{enumerate}}
\newcommand{\een}{\end{enumerate}}

\newcommand{\Z}{\mathbb Z}
\newcommand{\R}{\mathbb R}
\newcommand{\N}{\mathbb N}

\newcommand{\s}{\ensuremath{\mathcal{S}}}

\newcommand{\om}{\ensuremath{\omega}}
\newcommand{\Om}{\ensuremath{\Omega}}

\newcommand{\La}{\ensuremath{\Lambda}}
\newcommand{\si}{\ensuremath{\sigma}}

\begin{document}

\title{{\bf Markov and almost Markov properties \\in one, two or more directions.}}

\author{ Aernout C.D.  van Enter \footnote{ Bernoulli Institute, University of Groningen, Nijenborgh 9, 9747AG,Groningen, Netherlands,
 \newline
 email: aenter@phys.rug.nl},\\ 
  Arnaud Le Ny \footnote{Univ. Paris Est Creteil, CNRS, LAMA, F-94010 Cr\'eteil, France,
  \newline
 email:  arnaud.le-ny@u-pec.fr},\\
Fr\'ed\'eric  Paccaut \footnote{ LAMFA CNRS UMR 7352, Universit\'e de Picardie Jules Verne, 80039 Amiens, France,
\newline
email: frederic.paccaut@u-picardie.fr},
}

\maketitle

\begin{center}
{\bf Abstract:} 
\end{center}

In this review-type paper written at the occasion of the Oberwolfach workshop {\em One-sided vs. Two-sided stochastic processes} (february 22-29, 2020), we discuss and compare  Markov properties and generalisations thereof  in more directions, as well as weaker forms of conditional dependence, again either in one or more directions. In particular, we discuss in both contexts various extensions of Markov chains and Markov fields and their properties, such as $g$-measures, Variable Length Markov Chains, Variable Neighbourhood Markov Fields, Variable Neighbourhood (Parsimonious) Random Fields, and Generalized Gibbs Measures.

\medskip

\footnotesize



\normalsize
\section{Introduction}

The first occurrence of the Markovian terminology in Probability Theory arose at the beginnning of the twentieth century,   when the Russian mathematician A.A. Markov investigated  dependent random variables in order to extend classical probability laws beyond independence, like at first   Chebyshev's law of large numbers. In the period 1907-1911 he started to study it at a fundamental level  \cite{Markoff1907,Markoff1907f}, and afterwards in a  popular application, published in 1913,  in which he managed to detect some regularities in the dependencies between vowels and  consonants in  A. Pushkin's novel "Eugene Onegin" \cite{Markoff1913}. In these articles, Markov summarized  his intensive work to investigate such "dependencies", for which he developed a technique going beyond (independent) coin tosses, by introducing a {\em time-oriented, of limited scope, principle} -- nowadays known as the {\em Markov property} which drives Markov chains -- :  Successive letters are considered to have some relationship, but when a given (present) letter is known, the next letter is independent of all the earlier ones. Later, during  the whole past century, and probably even more during the last twenty years, this "Markov chains" concept has been used in various areas of theoretical and applied probability and statistics (to study such varied topics as   DNA sequences, coding theory, weather predictions, etc.); the concept was also   extended to  sequences with longer dependence on the past (like  Doeblin-Fortet's "{\em Cha\^ines \`a liaisons compl\`etes}" \cite{DoeFor} or Berbee's "{\em Chains with infinite connections}" \cite{Berbee}), to stochastic processes (Browian motion and beyond), and it was also extended to  spatial dependence through the notion  of "Markov fields", implicitly already widely studied in statistical physics through {\em e.g.} the famous Ising model of ferromagnetism \cite{Ising}. This can be summarized briefly as follows:

In a {\bf Markov Chain} the future is independent of the past, given the present.

In a {\bf Markov Field} the inside of a (finite) area is independent of the outside, given the border. Or, differently put, two areas, one of which is finite,  are independent, given the border between them.  

It is well-known\footnote{Concerning finite-state discrete-time stationary contexts, it has been known since almost 50 years that  homogeneous Markov Chains are equivalent to   translation-invariant Markov Fields, see for example \cite{HOG}, Ch 3. This result  was independently obtained by Brascamp \cite{Bra} and Spitzer \cite{Spi}.} that discrete-time finite-state Markov Chains, which are described by conditional probabilities for events in the future, and thus condition in one direction, 
 can also be described as one-dimensional Markov Fields, that is,  nearest-neighbor Gibbs measures for finite-spin models, described by  conditional probabilities in two directions, left and right. In such Markov fields the time interpretation of past and future is  being replaced by the space interpretation of an interior volume, surrounded by an exterior to the left and to the right.

If space is  $\Z$, one-dimensional, and time is discrete (thus also the one-dimensional  integer line $\Z$), the difference in description  is therefore between a conditioning in one direction (time, the past) versus a conditioning in two directions (space, left and right). But despite this, one obtains the same class of measures, as long as the conditioning is Markovian.

 On the other hand, if we consider arbitrary stochastic processes, one-directional  descriptions and two-directional descriptions can provide highly non-equivalent results.  
 In short, as observed before, in long-range settings, {\em neglecting  what lies outside borders is not equivalent to forget history, except for short-sighted approximations} (see {\em e.g.} \cite{vEB} for a recent review).
 
 If we relax the Markov requirement to weak dependence, that is, continuous --in the product topology-- dependence,  either on the past (generalising the Markov-Chain description) or on the external configuration (generalising the Markov-Field description), it turns out this equivalence breaks down, and neither class contains the other. In one direction this result has been known for a few years, in the opposite direction a counterexample was found more recently. The counterexample in this direction  is based on the phenomenon of entropic repulsion in long-range Ising (or "Dyson") models. In higher dimensions, even in the Markov case, there is no equivalence.  Moreover, there are different ways of considering Markov properties,  where one conditions on finite or infinite boundary configurations, on lexicographic pasts, etc.
 
 On trees,  an often-considered version of a one-directional, global,  Markov property provides one with the class  of "splitting" or "Markov chain" Gibbs measures \cite{HOG,Roz}.

Here we focus on the  diffusion of these extended notions within the communities at the interface of computer science (queueing networks), mathematics (probability and statistics, statistical mechanics, dynamical systems, ergodic theory) and theoretical physics, focusing on the mathematical point of view.  We discuss some of the known issues, properties and  counterexamples. More precisely, we investigate and compare their extensions in space and time in the one-sided {\em vs.} two-sided direction, {\em i.e.} we compare various extensions of Markov chains {\em vs.} Markov fields.

 We also discuss some further ways of weakening the continuity properties to almost sure continuity, emphasing systems  in one dimension, which for example play a role in describing systems with non-compact state spaces, systems with unbounded memories (within the VLMC framework described below), or coarse-grained models in phase transition regions ("RG pathologies" which giving rise to non-Gibbsianness and Generalized Gibbs measures within the so-called {\em Dobrushin program of restoration of Gibbsianness}).

We will thus concentrate on systems having  either Markov, or at least some continuity --in the product topology-- (also known as "almost Markov", cf Sullivan \cite{Sul}), properties.
The class of systems with continuity properties in one direction are the $g$-measures, those with continuity properties in two directions form  the class of Gibbs measures. Such continuity properties in space  in fact characterise {\em Gibbs measures}, also in higher dimensions \cite{Koz,Sul, BGMMT}. 
Measures with such continuity properties tend to be better behaved. For example, in Gibbs measures, as opposed to general stochastic processes,  \cite{Gur,OW,BuSt} the one-directional Kolmogorov-Sinai entropy and the two-directional entropy density coincide.\footnote{We remark that the quantity which we here call a ``two-directional`` entropy, is a one-dimensional  example of  entropy-like quantities which can be defined on more general graphs; such quantities have been  also called  ``inner'' or ``conditional'' or ``lower'' or "erasure" entropies in the literature \cite{AvELMP,EV, FS,Temp,VW}.}
 
In fact, not only  is it known that one-directional and two-directional entropy densities for Gibbs measures (for absolutely summable potentials)  are identical \cite{Temp, AvELMP}, but under some stronger conditions on the interaction decay, one-directional and two-directional tail properties of Gibbs measures are the same \cite{HolSt}.

In one dimension, it was open for a long time if the  $g$-measure property (one-directional) and the  Gibbs-measure property (two-directional) were equivalent or not.
 The measures which have continuous conditional probabilities in one direction are nowadays known by different names (even if  continuity was not a part of  most of the original definitions); not only  as $g$-measures, but also as ``chains with complete connections'' or `` chains of infinite order'', or ``random Markov chains''. They were introduced in the thirties and repeatedly rediscovered (under different names) \cite{DoeFor, Ruma, Harr,Keane, Kalik}. A few years ago, Fern\'andez,  Gallo and Maillard \cite{FGM} constructed a $g$-measure -- with one-directional continuous conditional probabilities-- which is not a Gibbs measure, as its two-directional  conditional probabilities are not continuous.

In  \cite{ BEEL}  an opposite result was found,  namely  that the Gibbs measures of the Dyson models --which, although of infinite-range type, have two-sided continuous conditional probabilities-- are not $g$-measures, as their one-sided conditional probabilities are not continuous. 
In contrast to most  counterexamples,  in this case the two-directional   behaviour is more ``regular'', more ``stochastic'',  than the one-directional behaviour. 

In the next section we describe some of the general theory and some examples and counterexamples in one dimension. In the following section we discuss some possible extensions of Markov and almost Markov properties to higher dimensions. After this we discuss the effects of weakening continuity properties even more by discussing non-Gibbsianness issues and reviewing the notions of {\em Variable Length Markov Chains} (VLMC), Parsimonious or {\em Variable Neighbourhood Random Fields} (VNRF), {\em Generalized Gibbs Measures} (GGM), {\em i.e.} other types of (non-uniform) longer ranges.

\section{Two directions in one dimension, background and notation}

\subsection{Specifications and Measures}

The classical notions of a Markov Chain and a Markov Field can be widely generalised  using topological (continuity) extensions of Markovian properties (weak dependence on the past in some topological sense), as follows.

 The generalization of Markov Fields leads to the Gibbs (or "thermodynamic") formalism. We refer to \cite{Bov, HOG,EFS,Fer,FV, Rue} for   general treatments of the Gibbs formalism. The Gibbs formalism was developed to provide a mathematically rigorous description of statistical mechanics, at the end of the 1960's. It also provides a rigorous framework to formalize properly the notion of  Equilibrium States, having  thermodynamic properties based upon  Boltzmann's original description and\textbackslash or  Gibbs's original ensembles \cite{Gibbs, Boltzmann}. The main model which has been studied in depth, the famous Ising model \cite{Ising}, is a prototype of Markov field.

 We consider  spin  models on lattices ${\Z}^d$, in the case $d=1$ for this section,  with Gibbs measures on infinite-volume product configuration spaces, $(\Omega,\mathcal{F},\rho)=(E^{{\Z}^d},\mathcal{E}^{{\otimes {\Z}^d}},\mu_o^{{\otimes {\Z}^d}})$. We mainly consider  translation-invariant models, occasionally on more general graphs than ${\Z}^d$.  We denote by 
$\mathcal{S}$ the set of  finite subsets of $\Z$ and, for any $\La \in \s$, 
write $(\Om_\La,\mathcal{F}_\La,\rho_\La)$ for the finite-volume configuration 
space $(E^\La,\mathcal{E}^{\otimes \La},\mu_o^{\otimes \La})$.  

  In the simplest case  the single-site state space is  the Ising space 
$E=\{-1,+1\}$,
with the a priori counting 
measure $\mu_0=\frac{1}{2} \delta_{-1} + \frac{1}{2} \delta_{+1}$.

Microscopic states or configurations, denoted by $\si,\om, \eta, \tau,\;$ etc., are  elements of  $\Omega$, 
equipped with the product topology of the discrete topology on $E$, for which these configurations are close when they coincide on large finite regions $\Lambda$ (the larger the region where they are equal, the closer the configurations are).

We denote by $C(\Om)$ the set of continuous (quasilocal) functions on $\Om$, which can be characterized here by uniform continuity:

\be \label{qlocfu} 
f \in C(\Omega) \; \Longleftrightarrow \; \lim_{\Lambda \uparrow \Z} \sup_{\sigma,\omega:\sigma_\Lambda=\omega_\Lambda} \mid f(\omega) - f(\sigma) \mid = 0.
\ee

In a number of our examples   we will consider {\em ferromagnetic pair} interactions, which  provides us with an extra tool: We can make use of {\bf FKG inequalities}.
Monotonicity for functions and measures concerns the natural partial (FKG) \cite{FKG} order "$\leq$ ", which we have on our Ising spin systems : $\sigma \leq \omega$ if and only if 
$\sigma_i \leq \omega_i$ for  all  $i \in \Z$. Its maximal and minimal elements 
are the configurations $+$ and $-$, and this order extends to functions: 
$f:\Omega \longrightarrow \mathbb{R}$ is called {\em monotone increasing}  when 
$\sigma \leq \omega$ implies $f(\sigma) \leq f(\omega)$. For measures,  we write $\mu \leq \nu$ if and only
if $\mu[f] \leq \nu[f]$ for all $f$ monotone increasing.
FKG arguments are based on the observation that increasing the interaction,  or increasing an external field, will FKG-increase the associated Gibbs measure. 

\smallskip
Macroscopic states are  represented by probability measures  
on  $(\Omega,\mathcal{F},\rho)$. At positive temperatures, 
following Dobrushin, Lanford and Ruelle, DLR (or  Gibbs) measures are defined in terms of consistent systems of (regular versions of) finite-volume conditional probabilities, for finite-volume configurations with prescribed boundary conditions outside of those volumes \cite{Dob1,LaR}.

Such a family of everywhere, rather than almost everywhere, defined
 conditional probabilities is called a {\em specification}.
  There is an advantage to this replacement of "almost everywhere" by "everywhere",  and even a need for it, as one has not yet a measure,  which could determine that something is  almost sure with respect to it.
 A measure for which a version of its  conditional probabilities provides  those of the specification is said to satisfy the DLR conditions for that specification.  The specifications of interest in the theory of lattice systems usually satisfy a finite-energy condition of non-nullness type. This says that no configuration in a local volume is excluded (has probability zero), uniformly in the boundary conditions. Moreover, the condition of continuity (or quasilocality) is required. This condition says that all conditional probabilities are continuous (quasilocal) functions of the boundary conditions.

A {\it measure} is said to be  quasilocal when it is specified by a quasilocal specification, and the main examples of them are the Gibbs measures as defined below, in terms of physical and thermodynamical notions ("interaction", "energy", "entropy").

\begin{remark}
In contrast to Kolmogorov's theorem, which says that a system of consistent marginal probabilities defines precisely one probability measure, for finite-spin specifications the number of measures satisying the DLR conditions for that specification can be either one or infinitely many. The latter situation sometimes is taken as the definition of a {\bf Phase Transition} (of first order). 
\end{remark}

A particularly important approach to  quasilocal measures consists thus in considering  the {\em Gibbs measures} with (formal) Hamiltonian $H$ defined  via a potential (or interaction) $\Phi$, a family $\Phi=(\Phi_A)_{A \in \s}$ of local functions $\Phi_A \in \mathcal{F}_A$.  The contributions of spins in finite sets $A$ to the total energy define the {\em finite-volume  Hamiltonians with free boundary conditions} :
\be \label{Ham}
\forall \Lambda \in \s,\; H_\Lambda(\omega)=\sum_{A \subset \Lambda} \Phi_A(\omega),\; \forall \omega \in \Omega.
\ee
To define Gibbs measures, we require for $\Phi$ that it is  {\em Uniformly Absolutely Summable} (UAS), {\em i.e.} 
that $\sum_{A \ni i} \sup_\omega |\Phi_A(\omega)| < \infty, \forall i \in \Z$.

One then can give sense to the  {\em Hamiltonian at volume $\Lambda \in \s$ with boundary condition $\omega$} defined for all $\sigma,\omega \in \Om$ as $$H_\Lambda^\Phi(\sigma | \omega) := \sum_{A \cap \Lambda \neq \emptyset} \Phi_A(\sigma_\Lambda \omega_{\Lambda^c}) (< \infty).$$

The {\em Gibbs specification at inverse temperature $\beta>0$} is then defined by
\be \label{Gibbspe}
\gamma_\Lambda^{\beta \Phi}(\sigma \mid \omega)=\frac{1}{Z^{\beta \Phi}_\Lambda(\omega)} \; e^{-\beta H_\Lambda^\Phi(\sigma | \omega)} (\rho_\Lambda\otimes \delta_{\omega_{\Lambda^c}}) (d \sigma)
\ee
where the partition function $Z_\Lambda^{\beta \Phi}(\omega)$ is a normalizing constant. Due to the
UAS condition, these specifications are quasilocal. It turns out that the converse is also true up to a non-nullness condition (finite energy) (see {\em e.g.} \cite{HOG, Fer, Koz, Sul, ALN2}) and one can take :
\begin{definition}[Gibbs measures]:\\
$\mu \in \mathcal{M}_1^+$ is a Gibbs measure iff $\mu \in \mathcal{G}(\gamma)$, that is, the conditional probabilities of $\mu$ - at least  a version thereof- are those given by $\gamma$  where $\gamma$ is a non-null and quasilocal specification, so that
\end{definition}

\be \label{esscont}
\lim_{\Delta \uparrow \mathbb{Z}} \sup_{\omega^1,\omega^2 \in \Omega}  \Big| \mu \big[f |\mathcal{F}_{\Lambda^c} \big](\omega_\Delta \omega^1_{\Delta^c}) - \mu \big[f |\mathcal{F}_{\Lambda^c} \big](\omega_\Delta\omega^2_{\Delta^c})\Big|=0
\ee
Thus, for Gibbs measures the conditional probabilities always have continuous versions, or, equivalently,
 there is no point of essential discontinuity. Points of essential discontinuity are  configurations which are points of discontinuity for ALL versions of the conditional probability. In particular one cannot make conditional probabilities continuous by redefining them on a measure-zero set if such points exist. In the generalized Gibbsian framework, one also says that such a configuration is a {\em bad configuration} for the considered measure, see e.g. \cite{ALN2}.
The existence of such bad configurations implies non-Gibbsianness of the associated measures, which have been mostly detected and identified starting from the late 1970's, in the context of Renormalization Group (RG) pathologies, which were studied in theoretical physics; RG theory was mainly developed for the analysis  of critical phenomena, for an early review of those pathologies see \cite{EFS}.

\begin{remark}
 If the interaction is of finite range (a Markov Field according to an extended definition), or sufficiently fast decaying, uniqueness of the Gibbs measure holds; indeed  no phase transition is expected in one dimension in considerable generality.
  But the Dyson models discussed below provide a counterexample of infinite range (they are polynomially decaying long-range Ising models), which however still satisfy the  UAS condition. 
\end{remark}

\begin{remark}
In fact it is enough to know the single-site conditional probabilities -- the single-site specification --, all other conditional probabilities can be obtained from those (especially in translation-invariant contexts). 
\end{remark}

Although the extension of the DLR equation to infinite sets is direct in case of uniqueness of the DLR measure for a given  specification \cite{FP, Foll,Gold2},  it can be more problematic otherwise: it is valid for finite sets only and  measurability problems might arise in case of phase transitions  when one wants to extend them to infinite sets. Nevertheless, 
 beyond the uniqueness case, such an extension was made possible by Fern\'andez and Pfister  \cite {FP} in the case of attractive models, that is models satisfying FKG properties. 
 
The concept they introduced is that of  a {\em global specification}, which appeared to be crucial to properly deduce one-sided conditional probabilities from two-sided (DLR) ones.

A  global specification is a set of consistent conditional probabilities where one considers probabilities of sets which have their supports not only in finite sets, but in more general sets $\Delta \in \Z$, which can be infinite, possibly  with infinite complements. The existence of such a global specification can be invoked to derive the existence of conditional probabilities of sets in  $E^\Delta$, and the possibility of conditioning  on 
configurations in $E^{\Delta^c}$. 

Note, by considering $\Delta=\mathbb{Z}$, that the set of measures a version of whose conditional probabilities is given by a global specification  contains at most one element. \\
The case we will be most interested in is the situation where we condition on only one half-line.
\noindent
This leads us to the concept of $g$-measures.

 The formalism of $g$-measures can be  developed in a parallel manner to the Gibbs formalism, but only using one-directional objects (conditional probabilities, specifications, etc.). \\
It was developed in probability and dynamical systems, and has its natural habitat in one dimension, as opposed to the more general Gibbs formalism.\\
We will call a measure a $g$-measure, once the future depends in a continuous manner on the past. Thus $g$-measures form a natural generalisation of Markov {\em Chains}, while Gibbs measures are a natural generalisation of Markov {\em fields}.

\begin{definition}
Let $\mu$ be a  measure on $\Omega ={E}^\Z$. We will call $\mu$ a $g$-measure for the function $g$ if the conditional probability for the next symbol being $a$, $\mu(x_0 =a| \{\omega_i \}_{i< 0}) = g(\omega_{{\Z}^{-}}a)$, depends in a continuous manner on the past, that is,  $g$ is continuous function on ${E}^{\Z^{-}}$.  

\end{definition}  
For a more extensive  description and background,  and also a comparison between $g$-measures and Gibbs measures in one dimension, we refer to \cite{BFV}.
 \begin{remark}
 We can also obtain the equivalent of a specification, a ``LIS'' (Left Interval Specification) for which the measure is a $g$-measure, analogously to the Gibbs measure definition \cite{FM1, FM2}. This is based on the observation that one can build general conditional probabilities from single-site ones.
\end{remark}

 {\bf Phase Transitions for $g$-measures}: 
 
 Although in the Markov-Chain set-up (which is a simple and well-known example of a $g$-measure) uniqueness holds, in the general $g$-measure case, similarly to the Gibbsian set-up, under only  the  condition of continuity phase transitions are possible \cite{BK,BHS}. However, although under appropriate uniqueness conditions (uniform boundedness of boundary energies,  or Dobrushin uniqueness, for example) $g$-measures and Gibbs measures turn out to be the same objects (\cite{BFV,FM1,FM2}, no general equivalence seems to apply in a more general setting. Thus the connection between the different classes of phase transitions possible is mostly unknown (if there  even is one). 

Proofs of phase transition, as well as of non-Gibbsianness, for $g$-measures,  so far seem to require a somehow more "constructive" approach to find the existence of (counter-)examples of the needed $g$-functions. They may be  employing a  condition of sparseness, or a condition of alternating sums  for sequences converging to zero, to determine the probabilities and the rules for   the distances one looks back into, or some other quite explicit construction. For example, in \cite{BK,BHS}, with a certain probability one may look at a certain distance in the past, to follow some majority-like rule on blocks; the chosen distances form a sparse sequence, and the probabilities of looking that far are carefully chosen in such a way that continuity still holds.  \\
The Gibbsian Dyson models which follow,  with their monotone, polynomial decay, have a more "natural" feel to them. The proofs of phase transitions for them, which are mentioned below, tend to have quite a different flavour.
A more "natural" example where one-directional conditional probabilities seem better behaved than two-directional ones, is the {\bf  Schonmann projection}, the marginal of an extremal  low-temperature 2-dimensional Ising model Gibbs measure to the one-dimensional line $\Z$. This is known not to be a Gibbs measure \cite{Sch,EFS}, due to a wetting phenomenon in the direction perpendicular to the line $\Z$. 
There are indications that it could be a $g$-measure \cite{BC}. Whether this is really so, is still open, however. Some of the used arguments came from a paper which already showed that one-directional conditional probabilities can behave in a different and more regular manner   than two-directional ones,  
 by Liggett and Steif  in  \cite{LiSt}, especially  Theorem 1.2 and Proposition 2.1. It appears to be the case, that, in contrast to the Dyson models, there is no wetting and entropic repulsion if a single interval of minuses is inserted in the plus phase, whereas entropic repulsion and wetting  do occur, via a merging of two wet regions which have increasing size in the perpendicular direction, once two large enough minus intervals are inserted. The distance of those intervals can be arbitrarily  large, as long as their lengths  are large enough (and much larger than the distance between the two intervals). 


\subsection {Long-range. Gibbs measures for Dyson Models and g-measures}

Here we describe some  properties of polynomial one-dimensional long-range Ising spin models, 
also known as {\em Dyson models}.

They provide the first example of a Gibbs measure which failed to be a $g$-measure, at low temperatures. However, at high temperatures, within the Dobrushin uniqueness region, they are known to be $g$-measures.  

In its original work, Dyson \cite{Dys} considered an Ising spin system in one dimension (on $\Z$), with formal Hamiltonian given by 
\begin{equation}\label{HamD}
H(\omega) = - \sum_{i>j} J(|i-j|)\omega_i\omega_j
\end{equation}
and $J(n) \geq 0$ for $n\in \mathbb{N}$ is of the form $J(n)=n^{-\alpha}$.

A conjecture due to  Kac and Thompson  \cite{Kac} had stated that there should be a phase transition for low enough temperatures if and only if $\alpha \in (1,2]$ (in zero magnetic field). Dyson proved a part of the Kac-Thompson conjecture, namely that for long-range models with interactions of the form $ J(n) = n^{-\alpha}$ with $\alpha \in (1,2)$, there is a phase transition at low temperatures.
Later different proofs were found, \cite{FILS,Joh,CFMP,ACCN} and also the case $\alpha=2$ was shown to have a transition \cite{FrSp}. 

In summary the following  holds (we call $\gamma^D$ the corresponding Gibbsian specification): 

\begin{proposition}
\label{DyFrSp} 
\cite{Dys,FrSp,Rue72,HOG,FILS,ACCN,CFMP,LP,Joh}.
The Dyson model with polynomially decaying potential, for $1< \alpha \leq 2$, exhibits a {\em phase transition at low temperature}: 
$$ \exists \beta_c^D >0, \; {\rm such \; that} \; \beta > \beta_c^D \; \Longrightarrow \; \mu^- \neq \mu^+ \; {\rm and} \; \mathcal{G}(\gamma^D)=[\mu^-,\mu^+] $$
 where the extremal  measures $\mu^+$ and $\mu^-$ are translation-invariant. They have in particular opposite magnetisations   $\mu^+[\sigma_{0}]=-\mu^-[\sigma_{0}]=M_0(\beta, \alpha)>0$ at low temperatures. Moreover, the Dyson  model in a non-zero homogeneous field $h$ has a unique Gibbs measure at all temperatures. 
\end{proposition}
It is well-known that  there is no phase transition for $J(n)$ being of finite range, and neither for $J(n)=n^{-\alpha}$ with $\alpha > 2$ (see {\em e.g.} \cite{BLP}). 
\begin{remark}
{\bf The borderline case:} \\The case of $\alpha=2$ is more complicated to analyse, and richer in its behaviour, than the other ones. There exists a hybrid transition (the "Thouless effect"), as the magnetisation is discontinuous while the energy density is continuous at the transition point. Moreover, there is a second transition below this transition temperature. In the intermediate phase there is a positive magnetisation with non-summable covariance, while at very low temperatures the covariance decays at the same rate as the interaction, which is summable. For these results, see  \cite{ACCN, I,IN}. 
\end{remark}
Note that another class of models, called {\em long-range Ising chains}, non translation-invariant and defined on the semi-infinite lattice $\N$, leads to other measures and phase transition regions,  with a borderline case at $\alpha =\frac{3}{2}$ \cite{FM1}.

\noindent
The proof in \cite{BEEL} uses  the approach of \cite{CFMP}, which has been extended to a number of other situations (Dyson models in random fields \cite{COP}, interfaces \cite{CMPR}, phase separation \cite{CMP17},  inhomogeneous decaying fields \cite{BEEKR}, etc). The disadvantage of this approach is that it works only at very low temperatures, as it is perturbative, and it works only for a reduced set of $\alpha$-values, $\alpha^* < \alpha <  2$,  with ${\alpha}^{*} = 3 - \frac{\ln 3}{ \ln 2}$. The advantage, however,  compared to other proofs,  is that translation-invariance does not play that much of a role. \\
The main idea of the approach of \cite {CFMP}, which was introduced in the $\alpha=2$ case by Fr\"ohlich and Spencer in \cite{FrSp}, is to construct a kind of triangular contours for which a Peierls-type contour argument can be obtained.  The energy of a contour of length $L$ has an energy cost associated to it of order $O(L^{2- \alpha})$, (and of order $O(\ln L)$ when $\alpha=2$).

\noindent 
Due to these  different behaviors in one dimension, there has been substantial interest in these Dyson models over the years. 
Varying the decay parameter $\alpha$ plays a similar role as varying the dimension in short-range models. This can be done in a continuous manner, so one obtains  analogues of well-defined models in continuously varying non-integer dimensions. This was one major reason why these models have attracted a lot of attention in the study of phase transitions and critical behaviour (see e.g. \cite{CFMP} and references therein).\\
For some  recent results  for these long-range Ising models with polynomially decaying interactions, see  \cite{LP, ELN, BEEKR, CMP17, EKRS, CELR,Ken}.

\smallskip
\begin{remark}

Some long-range $g$-measures  also  have been proven to have  phase transitions \cite{BK, BHS}, but those models and those proofs have a very different flavour from Dyson models. Although there even exist examples of $g$-measures with phase transitions which also are Gibbs measures \cite{BFV}, so far there seems to be no connection between the two situations.
\end{remark}
\medskip
The reason Dyson model Gibbs measures fail to be $g$-measures rests on the phenomenon of "entropic repulsion". 

 A main ingredient of this is  the occurrence of a mesoscopic interface localisation: \\
In \cite{CMPR} the main result describes the {\em mesoscopic} localisation of an ``interface point'' at low temperatures, for a Dyson model\footnote{satisfying extra technical restrictions on its parameters mentioned above, at low temperatures.} with Dobrushin boundary conditions (minus to the left, plus to the right).\\
 If one considers an interval $[-L,+L]$ with such Dobrushin boundary conditions, with probability close to 1 the interface point is near the center, that is, with large probability it is not more than 
$O(L^{\frac{\alpha}{2}})$
away from  the center with a (Gaussian) probability distribution; moreover the probability for the interface point  being at larger distances than $\varepsilon L$ from the origin  is bounded by $ O(L \exp -L^{2 - \alpha})$. The system is in a minus-like phase with strictly negative magnetisation left of the interface point, and in a plus-like phase with positive magnetisation to the right of the interface point.\\ 
From this penomenon it is intuitively plausible, and in fact not that hard to prove \cite{BEEL,vEB}, that the one-dimensional "Dobrushin boundaries" repel the interface point, due to entropic reasons. As such, it is a phenomenon which works at low temperatures, but not at $T=0$.  

Moreover for the case $\alpha =2$  there is  {\em no mesoscopic localisation}, but the position of the interface point has {\em macroscopic} fluctuations. Thus in that case our proof breaks down.

By a few small technical steps, it can be shown that the alternating configuration is a point of essential discontinuity  for the one-directional conditional probability.\\
In words,  this means that a large  alternating interval, preceded by a (VERY) large frozen minus interval configuration  again is followed by a minus phase, while, when it is preceded by a large frozen plus interval, it is followed by a plus-phase interval. But this dependence on the presence of  a frozen plus or minus interval far (of order $L_0$) to the left (= in the past), violates the continuity condition which is required for ${\mu}^{+}$ to be a  $g$-measure. As our  measure was defined to be a Gibbs  measure for the Dyson interaction, it automatically has two-sided continuity; thus we have obtained our counterexample.  

We can therefore conclude \cite{BEEL}:

\begin{theorem}

The very-low-temperature Gibbs measures of Dyson models decaying with a power $\alpha$  for $1 < \alpha < 2$ cannot be written as $g$-measures for a continuous $g$-function.  
\end{theorem}

Therefore, of the class of Gibbs measures for quasilocal specifications  and the class of $g$-measures with continuous $g$-functions, neither of the two classes contains the other one.

\section{Beyond one dimension: Continuity and  Markov properties in more directions}

\hspace{.5cm} {\bf Local versus Global Markov:}

\noindent
In higher dimensions, the analogue of the Markov-Field property occurs as a ``Local'' Markov property, where one conditions on configurations in the complements of finite sets.   
Other properties, such as the Global Markov property, which allows one to conditioning on boundary configurations of infinite sets with infinite complements like half-spaces, conditioning on  lexicographic pasts, or  being a (tree-indexed) Markov Chain (a ``splitting Gibbs measure'') on tree graphs, play a role which is  more like an analogue the Markov-Chain property. It is known, however, that in contrast to the one-dimensional situation  these properties do not follow from the Local Markov properties. 

Local Markov fields which violate Global Markov properties include mixtures of left-right and right-left Dobrushin states, mixtures of tree-indexed Markov chains, and there exist even translation-invariant extremal Gibbs measures without the Global Markov property. Most of the results on Global Markov properties mentioned in this section can be found in \cite{Foll, Gold2, GKS, Isr,Roz, HOG}. 

{\bf Global Markov and dependence on the lexicographic past:}

\noindent
 If one tries to compare such non-local Markov properties with nonlocal  continuity properties, by for example considering continuity properties as a function of the {\em lexicographic past}, (as a possible analogue of the one-dimensional ordinary past), one finds that they can be quite different. 
 
 Indeed,  it is not difficult to see that the low-temperature plus-phase of the Ising model in $d=2$, for example, although it is actually known to be well-behaved enough that it even satisfies the Global Markov Property, has conditional probabilities which display an essential point of discontinuity as a function of the lexicographic past. So even a Global Markov Field can have conditional probabilities which are discontinuous as a function of the (lexicographic) past. The proof of this statement is very close to that of the non-Gibbsianness of the Schonmann projection \cite{Sch}  (the marginal of the low-temperature two-dimensional Gibbs measure on the configurations a line $\{\Z,0\}$) as given in \cite{EFS}, and further studied in \cite{FP,BC} for example; only we now replace the line $\{\Z , 0 \}$ surrounding the origin by two half-lines $ \{{\Z}^-,0 \} $ and $ \{ {\Z}^+ \cup  0,-1 \}$ in the lexicographic past of the origin, then the proof goes through more or less literally. The wetting phenomenon which is responsible is identical: there is an entropic repulsion from a frozen interval into the ``future'' direction producing a wet droplet, and having two intervals left and right which are large enough, causes the two wet droplets to merge. Under conditions of strong uniqueness (high-temperature Dobrushin uniqueness, or Dobrushin-Shlosman conditions e.g.), continuity of the magnetisation in the origin as a function of the lexicographic past  configurations holds, however. (AvE thanks Brian Marcus and Siamak Taati for first asking him  about dependence properties on the lexicographic past and for discussions on this issue).  

{\bf Convex combinations need not respect Global Markov properties:}

Convex combinations of Gibbs measures which satisfy a Global Markov property, or are tree-indexed Markov chains, may, but need not, satisfy these properties.\\ Note that Global Specifications only allow for a single Gibbs measure satisfying them, and also that convex combinations of Gibbs measures for different Local (Gibbsian) Specifications never are Gibbs measures themselves, so they certainly don't satisfy a  Local Markov property. Thus a convex combination of two measures which satisfies a Global Markov property is specified by a Global Specification which differs from the two Global Specifications of those two measures.

{\bf Entropies:}

Analogues of the Kolmogorov-Sinai entropy in higher dimensions for Markov random fields, in which again conditioning on a measure in a restricted number of directions,  or on lexicographic pasts,  takes place, also have been considered, for example in \cite{ABMP,GKS}. 
Existence of such entropies turns out to be a  weaker property  than continuity.
\section{Beyond continuity}
\subsection{Weakening Gibbsian properties to full-measure sets: Continuity of conditional probabilities and  summability of interactions.}

 In Statistical Mechanics, it appeared a few decades ago that the above two-sided Gibbsian framework had to be relaxed a bit in order to incorporate renormalized measures and to get a more relevant framework to study scalings near criticality, as well as near first-order transitions, see  {\em e.g.} \cite{EFS}. Depending whether one focuses on a topological characterization (essential continuity of conditional probabilities) or on more thermodynamical ones (equilibrium states and uniformly absolutely summable  potentials), this gave rise to two main restoration notions, within a Dobrushin program launched in 1995 during a famous talk in Renkum (NL), namely {\em almost Gibbsianness} and {\em weak Gibbsianness}. Some of Dobrushin's ideas about these issues have appeared in  \cite{DS1,DS2}. While the latter notion (of weak Gibbsianness)  seemed indeed too weak to provide a proper notion of equilibrium states for renormalized measures, the former notion  has been  developed in the fundamental example of decimated Ising measures, with promising results, although it appears also still to be a bit too strong a notion. Intermediate  notions exists in this still active Dobrushin program of restoration of Gibbsianness (see {\em e.g.} the notion of intuitively weak Gibbs \cite{EV}, or also  the study of Tjur points to extend the region of Gibbsian measures to more parameters, as studied in  \cite{Ber}), and more studies are nowadays required for such examples and for the better understanding of {\em Generalized Gibbs measures}. 

In both  directions, the relaxation of the Gibbs property takes place in requiring the chosen  defining property (continuity of conditional probabilities or summability of interactions) to hold almost everywhere instead of for every configuration.

To settle the first notion, let us call  $\Omega_\gamma$ the set of {\em good configurations} of
a specification $\gamma$, which gathers its points of continuity.

\begin{definition}[Almost Quasilocal (or Almost Gibbs)]
A probability measure $\mu \in \mathcal{M}_1^+$ is said to be {\em almost Gibbs} if there is a specification $\gamma$ s.t. $\mu \in \mathcal{G}(\gamma)$ and $\mu(\Omega_\gamma)=1$.
\end{definition}

\begin{definition}[Weak Gibbs] A probability measure $\nu \in \mathcal{M}_1^+$ is  {\em Weakly Gibbs}
if there exists a potential $\Psi$ and a tail-measurable\footnote{Tail-measurability is required to insure that the partition
function is well-defined.} set 
$\Omega_\Psi$ on which  $\Psi$ is absolutely convergent, of full measure
 $\nu(\Omega_\Psi)=1$, s.t. $\Psi$ is consistent with $\nu$, i.e. $\mu \in \mathcal{G}(\gamma^{\beta \Psi})$ for some $\beta >0$.
\end{definition}

While most of the renormalized measures seem in fact to be weakly Gibbs, 
 this notion appears not quite strong 
enough to fully describe a satisfactory notion of equilibrium states without some extra topological requirements (\cite{KLNR}, see also \cite{LT}).   {\bf Almost Gibbsianness implies weak Gibbsianness} and this can be seen by a general  construction following the lines of the construction of  Kozlov's potential built from a vacuum potential (see \cite{Fer,ALN,MMR, MRVS}). The basic example showing renormalization group pathologies, decimation of Ising models in dimension two, has been shown to be almost Gibbs at low temperatures, but important counterexamples have been exhibited (and thus some work remains to be done to continue to fix a proper framework for equilibrium states, including renormalized measures).


\subsection{Weakening g-measure properties to full-measure sets: Continuity of g-functions.}

If $g$ is a continuous function on ${\Z}^-$, the existence of a $g$-measure compatible with $g$ is straightforward, using a fixed point theorem (for which compactness and continuity are required). Just as in the Gibbsian framework, weakening the continuity assumption requires the notion of essential discontinuity. Any measure $\mu$ on $\Omega$ defines a $g$-function $g_{\mu}$ by taking the limits $\mu(x_0=a|\{w_i\}_{-n\leqslant i<0})$ but $g_{\mu}$ is defined only $\mu$-almost surely. Likewise, if $g$ is given, being a $g$-measure $\mu_g$ for $g$ only involves a set of $\mu_g$ measure one. A discontinuity point of $g$ is \emph{essential}, if it cannot be removed by modifying $g$ on a set of measure zero. A function $g$ being given, it seems natural to think that the existence of a $g$-measure will hold if the set of discontinuity points is "small". The smallness of this set can be measured in several ways: by its shape (viewed as a set of infinite paths) as in \cite{GG} or by assuming its topological pressure to be strictly negative (it mixes the shape of the set with the values $g$ takes on this set) as in \cite{GP}. More generally, the existence is proved to hold in \cite{FGP} if the set of discontinuity points is being given zero measure for a dynamically defined candidate measure. But wilder things may happen. There exists a positive $g$-function (everywhere discontinuous) with no $g$-measure;  an example, due to Noam Berger, is given in \cite{FGP}. In the same paper is also exibited a $g$-measure which is proved to be everywhere essentially discontinuous; this example has first appeared in Harry Furstenberg's thesis and has been published in \cite{Fur}. It has also been proved not to be almost Gibbs in \cite{LMVV}, the specification being everywhere discontinuous.

Just as in the Gibbs framework, when the set of essentially discontinuous points has measure zero, the measure could be called \emph{almost-$g$} but more studies are required to know whether the existing example of non-$g$ measures arising from statistical mechanics could be incorporated in this class. In the case of the Gibbs measures of the Dyson model in the non-$g$ region of parameters, it is still an open question, although one knows that there are at least uncountably many discontinuity points \cite{Endo}.The notion of \emph{weak-$g$} does not seem to exist yet, but should similarly be investigated by e.g. extension of Walters's class of potential \cite{Walters}. For this, an interaction would have to be constructed from, say, the Left Interval Specification (LIS) defined in \cite{FM1, FM2}, in the same way Kozlov's potential is built in the Gibbs setting.

\subsection{Variable Length Markov chains, Parsimonious Random Fields and Generalized Gibbs Measures}
In dimension one, there is a particular class of $g$-measures, called \emph{Variable-Length Memory chains} or \emph{Variable-Length Memory chains} (VLMC), defined as follows.

Given a $g$-function $g$, for any $\omega\in\Omega$, set
\begin{align*}
\ell^g(\omega)&:=\inf\{k\geqslant1:g(\sigma_{\Z^-}\omega_{[-k,-1]}a)=g(\eta_{\Z^-}\omega_{[-k,-1]}a)\,,\,\,\forall a\in E,\sigma\in\Omega,\eta\in\Omega\}
\end{align*}
(with the convention that $\ell^g(\omega)=\infty$ if the set is empty). Call $\omega_{[-\ell^g(\omega),-1]}$ the context of $\omega$ (this terminology was first used by J. Rissanen in \cite{Ris} for information theory purposes). It is the smallest suffix of $\omega_{\Z^-}$ we need to get the distribution of the next symbol according to $g$. The set $\tau^g:=\cup_{\omega}\{\omega_{[-\ell^g(\omega),-1]}\}$ is sometimes called the {\em context tree}, because it can be pictorially represented as a rooted tree in which each path from the root to a leaf represents a context. If $\sup_{\omega}\ell^g(\omega)=k<\infty$, the chain is a Markov chain of order $k$, but when $\sup_{\omega}\ell^g(\omega)=\infty$, the chain is no longer Markov. The set $\tau^g$ is sometimes assumed to be countable, which amounts to saying that $g$ is piecewise constant on a countable number of cylinders. By using the combinatorial structure brought by the context tree, precise results on existence and uniqueness of $g$-measures (in some cases necessary and sufficient conditions) are obtained in \cite{G} and \cite{CCNPP}, without requiring global continuity of $g$.

\par This parsimonious description of one-dimensional chains using contexts has been extended to higher-dimensional fields in \cite{LO}: the \emph{Variable-Neighbourhood Random Fields} are defined by giving, for each site, the minimal neighbourhood needed to predict the symbol at the given site. As in the one-dimensional case, this minimal neighbourhood depends on the outside configuration. L\"ocherbach {\em et al.} \cite{LO} called them {\em Variable-neighborhood random fields} or {\em Parsimonious Markov random fields}, providing natural examples falling into the following Markovian framework:
\begin{definition}[Parsimonious {\em Markov} Fields \cite{LO}]\label{DefParsMark}
Let $\gamma$ be a specification. One says that  $\mu \in \mathcal{G}(\gamma)$ is a {\em Parsimonious Markov Field} if for any $\Lambda\in \mathcal{S}$ and $\mu$-a.e. $\omega$, there exists a {\em finite context} $C=C_\Lambda(\omega)$ $\subset \mathbb{Z}^d$ such that
$\gamma_\Lambda(\cdot \mid \omega_{\Lambda^c})=\gamma_\Lambda(\cdot \mid \omega_C)$
and for all $\tilde{C} \subset \mathbb{Z}^d$, if $\gamma_\Lambda(\cdot \mid \omega_{\Lambda^c})=\gamma_\Lambda(\cdot \mid \omega_{\tilde{C}}) $ then $C \subset \tilde{C}$.
\end{definition}

According to this definition, there might be a set of realizations of $\mu$-measure zero so that $|C_\Lambda(\omega)|=\infty$. This is e.g. the case of another example from \cite{CGL}, the so-called {\em Incompletely observed Markov random fields}. As in the latter, we extend Definition \ref{DefParsMark} allowing the context to be infinite, providing then a new\footnote{While Definition \ref{DefParsMark} is more a new description of particular, mostly Gibbs, class of random fields.} family of  possibly non-Gibbsian random fields.

\begin{definition}[Parsimonious {\em Random} Fields \cite{CGL}]\label{DefParsRand}
 One says that  $\mu \in \mathcal{G}(\gamma)$ is a {\em Parsimonious Random Field} if Definition \ref{DefParsMark} holds for non necessarily finite contexts $C=C_\Lambda(\omega)$  $\subset \mathbb{Z}^d$.
\end{definition}

In \cite{CGL},  the question of finiteness of the interaction neighborhood is related to the absence/presence of phase transition in the underlying Markov field, but in a more intricate way than one might suspect at  first sight. The question of  Gibbsianness or non-Gibbsianness of this model has not yet been  investigated, although we believe that it might be non-Gibbs in the phase transition region by a mechanism similar to the one described in \cite{EFS} for stochastic RG transformations. In \cite{ALN}, one of us described how the decimation of the Ising model could be incorporated in an even larger class of parsimonious measures.

So far, it appears that the relationship between one-sided VLMC and two-sided VNMF notions in this context have not been considered, but it might be of interest to investigate whether -possibly under some extra conditions- there are equivalences or implications.
\section{Conclusion, Final Remarks}

We have shown that on the one-dimensional lattice $\Z$, one-directional and two-directional properties, and in particular continuity of conditional probabilities, not only are not equivalent, but that, moreover, neither of the two continuity conditions  implies the other one.\\ 
In higher dimensions and on trees, different, non-equivalent, generalisations of the Markov chain property are possible.\\
Weakening continuity at  all points to continuity on full-measure sets opens up a whole new set of questions, many of which so far are unanswered.


{\bf Acknowledgments:} 

We thank the MFO to let us participate in miniworkshop at Oberwolfach on the issue of one-sided and two-sided properties of stochastic systems.  We thank the various participants as well  other colleagues and in particular our coauthors on these issues for the many discussions on these and related topics,  and for all that we learned from them. In particular we are grateful to: 
 Noam Berger, Stein Bethuelsen, Rodrigo Bissacot, Peggy C\'enac, Brigitte Chauvin, Diana Conache, Loren Coquille, Eric Endo, Roberto Fern\'andez, Sandro Gallo, Frank den Hollander, Bruno Kimura, Christof K\"ulske, Piet Lammers, Eva L\"ocherbach, Brian Marcus, Nicolas Pouyanne, Pierre Picco, Frank Redig, Wioletta Ruszel, Senya Shlosman,  Cristian Spitoni, Jeff Steif, Siamak Taati and Evgeny Verbitskiy.

\end{document}